\documentclass[11pt,a4paper]{amsart}
\usepackage[utf8]{inputenc}
\usepackage[T1]{fontenc}
\usepackage[english]{babel}
\usepackage{amsmath}
\usepackage{amsfonts}
\usepackage{amssymb}
\usepackage{graphicx}
\usepackage{lmodern}
\usepackage[ left=3.25cm,right=3.25cm,top=3.5cm,bottom=5cm]{geometry} %
\usepackage{amsthm}

\theoremstyle{definition}

\newtheorem{remark}{Remark}

\newtheorem*{example}{Example}

\newcommand{\norm}[1]{\left\lVert#1\right\rVert}
\newcommand{\M}{\mathcal M}
\renewcommand{\to}{\rightarrow}

\title[Anisotropic interaction avoiding collisions]{An anisotropic interaction model with collision avoidance} 

\author[C.~Totzeck]{Claudia Totzeck}
\address[C.~Totzeck]{University of Mannheim \\ 68159 Mannheim, Germany}
\email[]{totzeck@uni-mannheim.de}

\begin{document}

\begin{abstract}
	In this article an anisotropic interaction model avoiding collisions is proposed. Starting point is a general isotropic interacting particle system, as used for swarming or follower-leader dynamics. An anisotropy is induced by rotation of the force vector resulting from the interaction of two agents. In this way the anisotropy is leading to a smooth evasion behaviour. In fact, the proposed model generalizes the standard models, and compensates their drawback of not being able to avoid collisions. Moreover, the model allows for formal passage to the limit 'number of particles to infinity', leading to a mesoscopic description in the mean-field sense. Possible applications are autonomous traffic, swarming or pedestrian motion.  Here, we focus on the latter, as the model is validated numerically using two scenarios in pedestrian dynamics. The first one investigates the pattern formation in a channel, where two groups of pedestrians are walking in opposite directions. The second experiment considers a crossing with one group walking from left to right and the other one from bottom to top. The well-known pattern of lanes in the channel and travelling waves at the crossing can be reproduced with the help of this anisotropic model at both, the microscopic and the mesoscopic level. In addition, the 'right-before-left' and 'left-before-right' rule appear intrinsically for different anisotropy parameters.
\end{abstract}

\maketitle

\section{Introduction}
Modelling of pedestrian dynamics is a challenging problem that has gained a lot of attention in the last decades \cite{BellomoDegondTadmor,MFGames,DegondHierachy,BellomoCrowdDynamics,Helbing1,HelbingMolnar,Yanagisawa}. The applied mathematics community as well as  engineers developed a variety of models that focus on different features, e.g.~collision avoidance by defining side-stepping procedures \cite{MTWsidestepping}, by minimizing appropriate cost functionals \cite{BailoCarrilloDegond,Borzi} or by computing collision times to anticipate a future collision and adjust the trajectories of the pedestrians involved accordingly \cite{BailoCarrilloDegond,Christiani,DegondVision}. 

Other models are more interested in social aspects such as the influence of groups on evacuation dynamics \cite{AlbiInvisible,Helbing1}. Moreover, comparisons to real data and statistical model fitting, e.g.~in evacuation dynamics, helped in the optimization of sites and obstacles \cite{Bode2,Bode1,ChristianiPeri,GomesMTW,HelbingExperiments,HelbingPanic,HelbingEmpirical,PiccoliTosin,SiebenSchumannSeyfried}. Especially, in evacuation dynamics it is important to include the structure of the surrounding into the model. This can be done with the help of the Eikonal equation \cite{improvedHughes,MTWHuges,KlarTiwariMahato,Goatin}. In addition, the interaction of pedestrians and cars is an important field of study in the interest of security and traffic management \cite{BorscheMeurer1,BorscheMeurer2,Roundabouts}.

Up to the author's knowledge the community widely agrees that isotropic interaction models that are used for the modelling of swarms of birds, schools of fish, sheep-dog interactions or opinion dynamics \cite{AlbiPareschi,Sheep1,CBO2,CBO1,MotschTadmor} are inappropriate for a detailed and realistic description of pedestrian dynamics \cite{BailoCarrilloDegond}. This lack of details is unfortunate since these models allow for formulations on the microscopic, mesoscopic and macroscopic scale and even the limiting procedures are well understood, see \cite{AlbiPareschi,CarrilloReview,smoothed,bodyattitutecoordination,DegondMotsch,Golse,JabinMeanfield,Vicsek} for an overview. A starting point for such a hierarchy in pedestrians dynamics is given in \cite{ChistianiPiccoli,DegondHierachy} where the relation of kinetic and macroscopic models for pedestrian dynamics is discussed.

Here, we are going to adapt isotropic interaction models for crowd dynamic in two dimensions in order to build an appropriate model for such dynamics including collision avoidance. The main objective of this article is to reveal this interesting link between potential-based isotropic interaction schemes and anisotropic schemes allowing for pattern formation. In fact, the introduction of only one additional parameter allows for this significant change in the dynamic.

First, the model is described on the particle level, then a formal limit to a kinetic equation is discussed. The main new key ingredient is the introduction of an anisotropy in form of a rotation applied to the force vector of the interaction. This rotation does not model a property, like a view cone, explicitly. One may think of it as the perception that adjusting the velocity helps to avoid a collision. This is a perception everyone undergoes unconsciously in his/her every day life, for example, when strolling around the city. In fact, the rotation allows us to model the idea of evasion in an elementary and smooth way. Moreover, it is straight forward to include interaction with obstacles, which can be represented at fixed positions with artificial velocities. The approach is inspired by \cite{Kreusser2,Kreusser3,Kreusser1}, where the formation of fingerprint patterns is studied with the help of an anisotropic force field. Further studies may consider the anisotropy parameter to be density or obstacle dependent. In the sense that the rotation becomes larger or smaller in locations where obstacles or high densities are present.

Beside the application to pedestrian dynamics, the approach is very general and can be used for many other interaction schemes as well. For example applications in swarms of autonomous agents, as discussed in \cite{TaylorLuzziNowzari}, are very interesting and important for traffic and transport in the future. 

A driving question in this research is whether the formation of lanes in a channel and travelling waves at a crossing as found in \cite{Appert-Rolland,Ranetbauer1,MTWsidestepping,Ranetbauer2}, can be reproduced with the help of the rotation anisotropy. The answer is affirmative, which the simulation results at the particle and the mean-field level will show.

Many models for pedestrian dynamics work with view cones or other sensorial perceptions to approximate collision times. Hence, they influence the future trajectories of the pedestrians to make the simulations more realistic and to avoid collision\cite{BailoCarrilloDegond,DegondHierachy,EversFetecauRyzhik, EversFetecauSun}. These approaches lead to anisotropies as well. 
Note that this kind of anisotropy has a different effect than the anisotropy discussed in the present paper. Indeed, considering view cones leads to an anisotropy in the sense that binary interactions are not symmetric. For example, one agent sees the other, and interacts with him or her, while the other one does not see the first and thus does not react. 
In contrast, the anisotropy introduced in the following rotates the vectors of the interaction forces for each of the two interacting pedestrians individually. Hence, if one pedestrian interacts with another, the second will also interact with the first and the corresponding force vectors are rotated symmetrically. 

An advantage of the anisotropic interaction induced by this rotation is that the dynamic can be written as system of ordinary differential equations (ODE). Further, the pedestrians are interacting with all their neighbours at the same time. In particular, there is no need to estimate times of collision and to choose a specific pedestrian to interact with first.  
Moreover, a formal mean-field limit leads to a formulation on the mesoscopic scale (PDE). Numerical results on the micro and mean-field level show the influence of the anisotropy for pattern formation.

The article is organized as follows: a detailed motivation and description of pairwise interactions is given in the next section. Based on the pairwise interactions, the microscopic model for the pedestrian dynamic is derived in Section~\ref{sec:micro}. The formal derivation of the mean-field limit is given in Section~\ref{sec:scaling}. As we are interested in realistic simulations, we discuss boundary conditions and handling of obstacles in Section~\ref{sec:boundaryObstacles}. Numerical schemes and their implementation are discussed in Section~\ref{sec:numerics}, before we discuss numerical results in Section~\ref{sec:simulations}. We conclude and give an outlook to future work in Section~\ref{sec:conclusionOutlook}.

	\section{Collision avoidance by anisotropy}
	We recall an isotropic interaction model to clarify the starting point of the modifications. Following Newton's Second Law, isotropic interaction models are second order ODE systems \cite{AlbiPareschi,Sheep1,smoothed,MotschTadmor,Vicsek}, given by
	\begin{align*}
		\frac{d}{dt} x_i &= v_i, \quad i=1,\dots, N, \\
		\frac{d}{dt} v_i &= -\frac{1}{N}\sum_{j\ne i} K(x_i,x_j,v_i,v_j) ,\qquad i=1,\dots, N, 
	\end{align*}
	where the force $K(x_i,x_j,v_i,v_j)$ can arise as gradient of a potential, for example the Morse potential \cite{Morse},  
		\begin{equation*}
			K(x_i,x_j,v_i,v_j) = \nabla_{x_i} R e^{-|x_i - x_j| / r} - A e^{-|x_i - x_j| / a}, \quad A,R,a,r > 0,
		\end{equation*}
		or an alignment force, as in the case of the Cucker-Smale model \cite{smoothed}
		\begin{equation*}
			K(x_i,x_j,v_i,v_j) = \frac{1}{(1 + \| x_i - x_j\|^2)^\beta} (v_j - v_i), \quad \beta \ge 0.
		\end{equation*}
	For reasons of well-posedness the system is supplemented with the initial data
	$x(0) = x_0, \; v(0) = v_0.$
	Here and in the following we use the notation $x = (x_i)_{i=1,\dots,N}$ with $x_i(t) = (x_i^1(t), x_i^2(t)) \in \mathbb{R}^{2}$ and $v= (v_i)_{i=1,\dots,N}$ with $v_i(t) = (v_i^1(t), v_i^2(t)) \in \mathbb{R}^{2}$ for all times $t$ and $i=1,\dots,N.$
	These systems of ODEs are used to model behaviour of swarms of birds, flocks of sheep or fish. 
	
	As we are interested in the patters arising in pedestrian dynamics. We first add some intelligence in form of a desired journey to the system. In fact, we add an acceleration term $u(x,v)$ which may depend on the current position or velocity of the agents. This leads to the system equations given by
	\begin{align*}
		\frac{d}{dt} x_i &= v_i, \quad i=1,\dots, N, \\
		\frac{d}{dt} v_i &= u(x_i,v_i) -\frac{1}{N}\sum_{j\ne i} K(x_i,x_j,v_i,v_j) ,\qquad i=1,\dots, N, 
	\end{align*}
	With the help of the acceleration term $u$ we can model the scenario of two particles approaching each other as it is shown in for $\lambda = 0$ in Figure~\ref{fig:influence of lambda}. The particles are pushing each other away from their desired velocity, leading to undesired behaviour.
	
	Now, we modify the model to allow for collision avoidance such that each particle can return to the desired velocity after the interaction. First, we discuss pairwise interactions before we extend the idea to the microscopic model for $N \in \mathbb{N}$ agents. From the standard attraction-repulsion-alignment schemes, see e.g.~\cite{smoothed,Vicsek,MotschTadmor}, we borrow the idea that the interaction of two agents is based on their distance or their velocities. Thinking of collision avoidance, a drawback of these models is that the force is aligned with the vector connecting the positions or velocities of the two agents, respectively. Indeed, it may happen that two approaching agents stop in front of each other and are not able to proceed with their journey.
In the following we generalize the well-known ODE based systems for particle interaction in order to derive a reasonable model allowing for collision avoidance. The main ingredient is a rotation matrix which is based on the current state information of each agent, such that no anticipation of states or collisions in the future are necessary. Naturally, each agent is interacting with all other pedestrians. 
Hence, using this anisotropic model allows for anticipation in a very elementary way. The details of the rotation used for the anisotropy are the following.

We let two agents, represented by positions, $x$ and $y\in \mathbb{R}^2$, and velocities, $v$ and $w \in \mathbb{R}^2$, adjust their direction of movement if their velocity vectors point towards each other. We recall that the angle between two vectors $v, w$ is given by $\arccos(\frac{v \cdot w}{\norm{v}\norm{w}})$. Hence, for two agents approaching each other, this  quantity is $\arccos(-1) = \pi.$ To be able to model the tendency to move to the left or the right, we introduce an anisotropy parameter $\lambda \in [-1, 1]$ 
	leading to the interaction angle 
\begin{equation}\label{eq:alpha}
	\alpha = \begin{cases} \lambda \arccos(\frac{v \cdot w}{\norm{v}\, \norm{w}}) , &\text{for } v\ne 0,w \ne 0 \\ 0, & \text{else} \end{cases}. 
\end{equation}
Note that $\alpha$ is well-defined, since the Cauchy-Schwartz inequality ensures $|v \cdot w| \le \norm{v}\,\norm{w}.$ 

\begin{figure}[htp]
	\begin{center}
	\includegraphics[scale=0.6]{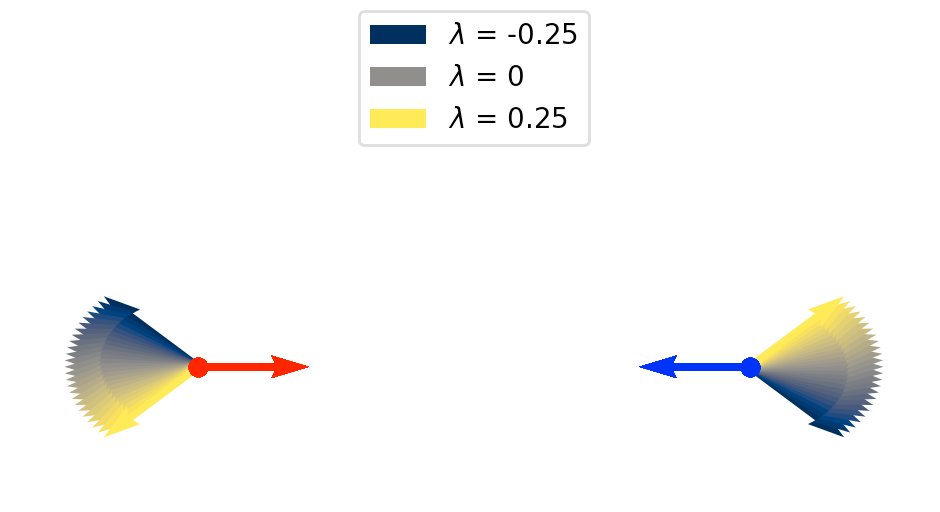}
	\caption{Visualization of the influence of $\lambda.$ For $\lambda = 0$ we are in the case of the grey vector. The interaction force reduces the force resulting for the desired velocity. The particles may stop in front of each other. For $\lambda > 0$ we are in the yellow region. The force vector resulting from the interaction is turned. Adding this rotated vector and the vector resulting from the desired velocity, yields the evasive behaviour. The red particle moves to the bottom and the blue particles moves to the top. For $\lambda < 0$ the roles of the to particles change. We are then in the blue region and the red particle moves to the top and the blue particle to the bottom.}
	\label{fig:visLambda}
	\end{center}
\end{figure}

\begin{figure}[htp]
	\begin{center}
	\includegraphics[scale=0.55]{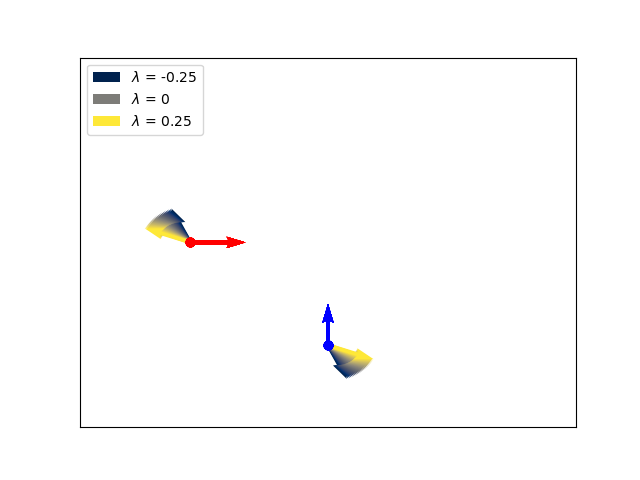}
	\includegraphics[scale=0.55]{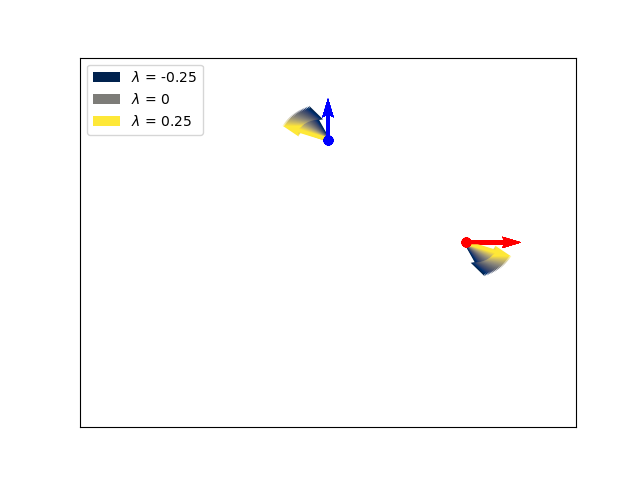}
	\caption{Visualization of the influence of $\lambda$ in a crossing scenario. Top: For $\lambda = 0$ we are in the case of the grey vector. The particles push each other diagonally and leave the path given by the desired velocity, see Figure~\ref{fig:influence of lambda} in the middle. For $\lambda > 0$ we are in the yellow region. The force vector resulting from the interaction is turned. Adding this rotated vector and the vector resulting from the desired velocity, yields the evasive behaviour. The red particle slows down and moves to the top before accelerating and adjusting the velocity to get back onto its path and the blue particle does not slow down as much as the red one and moves to the right. For $\lambda < 0$ the roles of the to particles change. We are then in the blue region and the red particle moves to the top and the blue particle moves to the right before being push back onto its path. Intrinsically, we see here the 'right-before-left' rule for $\lambda > 0$ and 'left-before-right' for $\lambda <0.$ }
	\label{fig:visLambda_Cross}
	\end{center}
\end{figure}

The resulting angle is now introduced into the rotation matrix that is multiplied to the force vector. In more detail, we have
\begin{equation}\label{eq:rotationmatrix}
	\M(v,w) = \begin{pmatrix} \cos \alpha & -\sin\alpha \\ \sin\alpha & \cos\alpha \end{pmatrix}.
\end{equation}
Note, that for $\lambda= 0$ the rotation matrix equals the identity. Hence, the model proposed here is a generalization of standard isotropic interaction models, see e.g.~\cite{smoothed,Vicsek}. 
	In the case of two agents, the equations read
	\begin{align*}
		\frac{d}{dt} x_1 &= v_1, &\frac{d}{dt} v_1 = u(x_1,v_1) - \frac{1}{2} \M(v_1,v_2) K(x_1,x_2,v_1,v_2),  \\
		\frac{d}{dt} x_2 &= v_2, &\frac{d}{dt} v_2 = u(x_2,v_2) - \frac{1}{2} \M(v_2,v_1) K(x_2,x_1,v_2,v_1).
	\end{align*}
In particular, approaching agents slightly turn the direction of the velocity vector and back-stepping may only happen, if the total force of approaching agents is very strong. More details on the influence of the parameter $\lambda$, in particular, the natural appearance of the 'right-before-left'-rule and the 'left-before-right'-rule, is visualized in Figure~\ref{fig:visLambda} and Figure~\ref{fig:visLambda_Cross}. Also, for agents walking in the same direction the matrix equals to the identity, leading to an isotropic interaction. As we assume to have small interaction ranges, the isotropic interaction takes place only when agents enter the comfort range of others. More details on this are given in the numerics section.

Altogether, the force acting on two agents positioned at $x$ and $y$ with velocities $v$ and $w$, respectively, is now given by the force coming for the standard interaction models multiplied with the rotation matrix.

\subsection*{Choice of $\lambda$}
	Clearly, the sign of $\lambda$ decides whether the agents have the tendency to move to the right ($\lambda > 0$) or to the left ($\lambda < 0$) to avoid a collision. The exact value of $\lambda$ is a modelling choice. Assuming that both agents reduce their speed w.r.t. the desired direction when approaching a possible collision, leads to  $\lambda \in (-1/2, 1/2).$ In fact, for $\lambda = \pm 0.5$ only one agents slows down, while the other keeps the speed into the desired direction but moves to the side.
	
	Assuming that the speed is changed significantly in a collision scenario motivates to choose $\lambda \in [-1/4, 1/4]$ as shown in Figure~\ref{fig:visLambda} and Figure~\ref{fig:visLambda_Cross}. In case of a crossing scenario $\lambda = \pm 1/4$ happens to imply naturally the 'right-before-left' rule or  'left-before-right' rule, respectively. This is the justification to choose $\lambda = 1/4$ for all numerical examples in the remainder. The anisotropic binary interaction is illustrated in the following example, before we state the microscopic model for arbitrary number of agents in the next section.

\begin{example}
	We consider two pedestrians in a channel and at a crossing. Therefore, let $u(x,v) = u_{r/b} - v$. Here, the indices $r$ and $b$ refer to a red and a blue pedestrian, and the desired velocities $u_r$ and $u_b$ will be fixed for the whole simulation. As interaction potential we choose the Morse potential as proposed in \cite{Morse}. 
	The influence of the rotation is illustrated in Figure~\ref{fig:influence of lambda}. In the first plot, the blue agent moves from the right to the left with $u_b = (-1,0)^T$ and the red agent from the left to the right with $u_r = (1,0)^T.$ The parameter $\lambda = -0.25, 0 , 0.25$ is varied. In Figure~\ref{fig:influence of lambda} (bottom) the crossing scenario is shown, the red pedestrians walks to the right, the blue one is walking from bottom to top with $u_b = (0,1)^T.$
	\begin{figure}[htp]
		\begin{center}
		\includegraphics[scale=0.5]{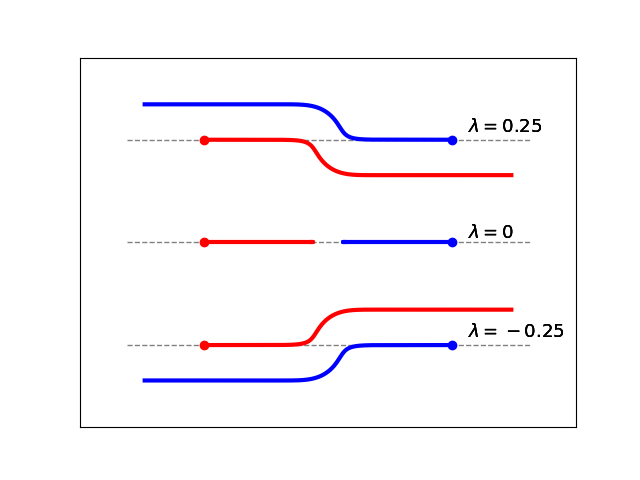}
		\includegraphics[scale=0.5]{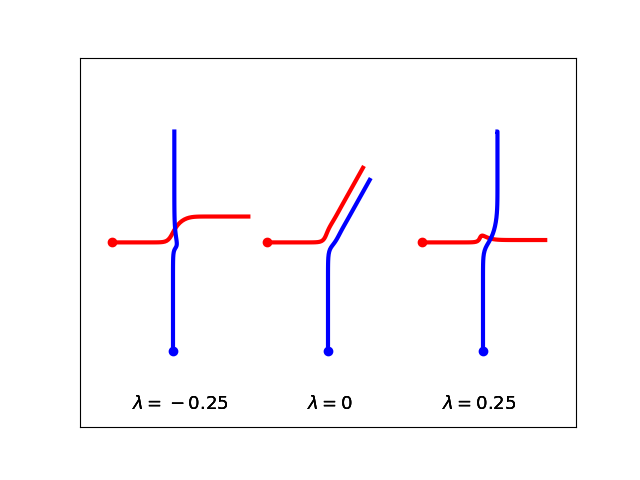}
		\caption{Illustration of isotropic and anisotropic pairwise interaction induced by different values of $\lambda$. The initial positions of the pedestrians are highlighted with a point. Top: Two pedestrians are approaching each other. Positive $\lambda$ leads to right-stepping, negative $\lambda$ results in left-stepping. The isotropic case corresponds to $\lambda = 0.$ Bottom: Two pedestrians meet at a crossroad. In the isotropic case, both depart from their desired trajectory. For $\lambda \ne 0$ the two return to their desired velocity after avoiding the collision. In case of the crossing with $\lambda = -0.25$ the blue pedestrians slows down, while the red pedestrian accelerates to avoid the collision. Then, the pedestrians relax their velocities towards the desired velocity. For $\lambda = 0.25$ the two pedestrians change their roles.
		}
		\label{fig:influence of lambda}
		\end{center}
	\end{figure}
	For $\lambda = -0.25$ the agents turn to the left, which is appropriate for countries with left-hand traffic. Similarly, for $\lambda = 0.25$ we have the case of right-hand traffic. For $\lambda = 0$ we are in the isotropic setting. We see that approaching agents stop in front of each other and are not able to pass. In Figure~\ref{fig:influence of lambda} (bottom) we see the scenario of two pedestrians with trajectories crossing each other. Again we observe the influence of $\lambda$. In the isotropic case the agents are not able to pass each other, in both other cases the pedestrians adjust their route and walk their way after a slight interaction. In case of the crossing with $\lambda = -0.25$ the blue pedestrians slows down, while the red pedestrian accelerates to avoid the collision. Then, the pedestrians relax their velocities towards the desired velocity. For $\lambda = 0.25$ the two pedestrians change their roles. Note, that the collision is avoided in all anisotropic cases.
	\end{example}
	
	\section{Microscopic Model}\label{sec:micro}
	In this section we introduce a general microscopic model for $N \in \mathbb{N}$ agents based on the pairwise-interactions discussed above.
	We introduce the anisotropy by rotating the force vector of the particle system derived from Newton's Second Law.  In particular, we obtain a second order ODE for each agent which reads
	\begin{subequations}\label{eq:ODE}
		\begin{align}
			\frac{d}{dt} x_i &= v_i, \quad i=1,\dots, N, \\
			\frac{d}{dt} v_i &= u(x_i,v_i) -\frac{1}{N}\sum_{j\ne i}   \M(v_i,v_j)\,K(x_i,x_j,v_i,v_j) ,\qquad i=1,\dots, N, \label{eq:velocityODE}
		\end{align}
		with anisotropy governed by the rotation matrix
		\begin{equation}
			\M(v_i,v_j) = \begin{pmatrix} \cos\alpha_{ij} & -\sin \alpha_{ij} \\ \sin\alpha_{ij} & \cos\alpha_{ij} \end{pmatrix}, 
		\end{equation}
		where the rotation angle is given by
		\begin{equation}
		\alpha_{ij} = \begin{cases} \lambda \arccos(\frac{v_i \cdot v_j}{\norm{v_i}\, \norm{v_j}}) , &\text{for } v_i\ne 0,v_j \ne 0 \\ 0, & \text{else} \end{cases},\quad 
		\end{equation}
		for $\lambda \in [-0.25,0.25]$ fixed. The system is supplemented with the initial data
		$x(0) = x_0, \; v(0) = v_0.$
	\end{subequations}
	
	Additionally, the model contains desired velocities $u(x_i,v_i)$ of each particle. The strength of interaction is described by $K$ and $\mathcal M$ rotates the vector in which the force acts on the agent. At this point we do not specify any boundary behaviour as it will strongly depend on the geometry under consideration. In the numerical section, we shall consider a channel and a crossing with periodic and reflecting boundary conditions. 
	
	Note that each agent interacts with \textit{all} other pedestrians. This is in contrast to many other collision avoiding models, where the interaction is considering only the agent $j$ which will collide with agent $i$ first. Moreover, the interaction depends on the position and the velocities of the agents, this increases the computational costs. We overcome this issue by limiting the interaction to a finite neighbourhood as is described in the following remark and in the section on the numerical schemes.
	
	\begin{remark}
		Note, that even though every agent is interacting with all other agents in this model. In reality, the interaction range is quite short. The range parameter in the interaction potential accounts for this behaviour. For the numerical simulations, we use the Morse Potential with short range interactions. Of course, the approach also works for models where agents interact only locally in a finite neighbourhood. Since the interaction depends on the positions and velocities of the agents, the model is computationally very costly. This is the reason why we use a finite range of interaction for the numerical simulations later on. See below for more details.
	\end{remark}	
	
	\begin{example}
		For the pedestrian scenarios in the numerics section, the above system can be specified as follows. We consider two groups of pedestrians, red ($r$) and blue ($b$), with $N_r$ and $N_b$ group members, respectively. We set $N:= N_r + N_b$. Each group has a desired velocity which we denote by $u_r$ for the red and $u_b$ for the blue. We use $u(x_i,v_i) = u_b - v_i$ for the pedestrians of the blue group and $u(x_i,v_i) = u_r - v_i$ for the pedestrians of the red group. Note that in this case the desired velocity depends only on the group each pedestrians belongs to and the current velocity of the agent. Moreover, we consider forces arising from Morse potentials that depend only on the distance of two interaction agents. This leads to $K(x_i,x_j,v_i,v_j) = K(x_i,x_j).$ Altogether, \eqref{eq:velocityODE} can be rewritten as
		\begin{subequations}\label{eq:ODE2groups}
			\begin{align}
				\frac{d}{dt} v_i &= (u_r-v_i) -\frac{1}{N}\sum_{j\ne i}   \M(v_i,v_j)\,K(x_i,x_j) ,\qquad i=1,\dots, N_r, \\
				\frac{d}{dt} v_i &= (u_b-v_i) -\frac{1}{N}\sum_{j\ne i}   \M(v_i,v_j)\,K(x_i,x_j) ,\qquad i=N_r+1,\dots, N.
			\end{align}
		\end{subequations}
	\end{example}
	\begin{remark}
		For $N \gg 1$ there is the well-known mean-field approximation of the interacting particle system in terms of measures which describe the probability of finding particles at a given location $x$ with velocity $v$. The mean-field approximation is formally derived in the following section.
	\end{remark}
	\begin{remark}
		We do not need a relaxation parameter in $u(x,v)$, as we can balance the terms with the help of the strength parameter of the interaction.
	\end{remark}
	
	\section{Formal derivation of mesoscopic model}\label{sec:scaling}
	In many interesting applications, like swarming or evacuation dynamics many agents are involved.  This motivates to seek for a mesoscopic description of the model which we derive in the following. Based on the agents' information given by \eqref{eq:ODE} we define the empirical measure
	$$f^N(t,x,v) = \frac{1}{N}\sum_{i=1}^N \delta(x_i(t) - x, v_i(t) -v).$$
	Let $\varphi$ be an arbitrary test function, we compute for $i=1,\dots,N$
	\begin{align*} 
		\frac{d}{dt} &\Big\langle f^N(t,x,v), \varphi(x,v) \Big\rangle = \Big\langle \partial_t f^N, \varphi \Big\rangle \\
		&= \frac1{N} \left[\; \sum_{i=1}^{N} \nabla_x \varphi(x_i(t),v_i(t)) \cdot v_i \right. \\ 
		&\left.\quad+ \nabla_v \varphi(x_i(t),v_i(t)) \cdot \Big( (u(x_i,v_i)-v_i) -\frac{1}{N}\sum_{j\ne i}   \M(v_i,v_j)\,K(x_i,x_j,v_i,v_j) \Big) \right] \\
		&= \frac1{N} \left[\; \sum_{i=1}^{N} \nabla_x \varphi(x_i(t),v_i(t)) \cdot v_i  \right. \\  
		&\left.\quad+ \nabla_v \varphi(x_i(t),v_i(t)) \cdot \Big( (u(x_i,v_i)-v_i) -\int  \M(v_i,\bar v)\,K(x_i,\bar x,v_i,\bar v)\, df_t^N(\bar x,\bar v) \Big) \right]\\
		&= \Big\langle \nabla_x \varphi(x,v) \cdot v  \\
		&\quad\qquad + \nabla_v \varphi(x,v) \cdot \Big( (u(x,v)- v) -\int   \M(v,\bar v)\,K(x,\bar x,v,\bar v)\, df_t^N(\bar x,\bar v)\Big)  ,f_r^N \Big\rangle. \\
	\end{align*}
	Now, assuming the existence of a limiting measure satisfying $f =  \lim\limits_{N\to\infty} f^{N} $ we formally pass to the limit $N \to\infty$ and use the variational lemma to obtain the evolution equation
	\begin{subequations}\label{eq:PDE}
		\begin{equation}
			\partial_t f + v\cdot \nabla_x f = \nabla_v \cdot \left[\left( \int  \M(v,\bar v) K(x,\bar x,v,\bar v)\, df(\bar x, \bar v) - (u(x,v)-v)  \right  ) f \right]
		\end{equation}
		supplemented with initial condition
		\begin{equation}
			f(0,x,v) = f_0(x,v).
		\end{equation} 
	\end{subequations}
	\begin{remark}\label{rem:meanfield}
		We emphasize that the mean-field equation cannot be derived in the case every particle has its own desired velocity as the particles would no longer be indistinguishable. One approach to prescribe a velocity field $u(x,v)$ is to use the solution of an appropriate Eikonal equation, see for example \cite{improvedHughes,MTWHuges,KlarTiwariMahato}. An investigation of a system coupled to an Eikonal equation planned for future work. Here, we adhere with the pedestrian scenarios discussed above. The details of the PDE in this specific case are discussed in the example below.
	\end{remark}	
	\begin{example}
		Under the assumption, that the crowd splits into two groups $f_r$ (red) and $f_b$ (blue) of $N_r$ and $N_b$ agents, respectively, we prescribe different desired velocities $u_r$ and $u_b.$ Without loss of generality we arrange the agents of the red group at the first $N_r$ positions of $x$ and $v$ and the blue group at the indices $N_R +1$ to $N.$ This leads to \begin{align*}
			u(x_i,v_i) &= u_r - v_i,\quad &&i=1,\dots,N_r, \\
			u(x_i,v_i) &= u_b -v_i,\quad &&i = N_r+1,\dots, N.
		\end{align*}
		Further, it holds
		\begin{align*}
			f^N(t,x,v) &= f^N_r(t,x,v) + f^N_b(t,x,v) \\&= \frac{1}{N} \left( \sum_{i=1}^{N_r} \delta(x_i(t) - x, v_i(t) -v) +\sum_{i=N_r + 1}^N \delta(x_i(t) - x, v_i(t) -v) \right).
		\end{align*}
		We observe that the convolution is linear w.r.t. $f^N$ to obtain the coupled PDE system given by
		\begin{subequations}\label{eq:pde-system}
			\begin{align}
				\partial_t f_r + v\cdot \nabla_x f_r &= \nabla_v \cdot \left[\left( \int  \M(v,\bar v) K(x,\bar x,v,\bar v)\, df(\bar x, \bar v) - (u_r-v)  \right  ) f_r \right]  , \\ 
				\partial_t f_b + v\cdot \nabla_x f_b &= \nabla_v \cdot \left[\left( \int \M(v,\bar v) K(x,\bar x, v,\bar v)\, df(\bar x, \bar v) -(u_b -v) \right) f_b \right].
			\end{align}
		\end{subequations}
		Note that the equations are coupled, as $f = f_r + f_b$ appears in the interaction terms.
		For the discussion of the numerical schemes later on, we introduce the notation
		$$ S_i(f) f_i := \left( \int  \M(v,\bar v)K(x,\bar x,v,\bar v)\, df(\bar x, \bar v) - (u_i-v)  \right) f_i, \qquad i \in \{r,b\}.  $$
		Then the evolution equations reduce to
		\begin{subequations}\label{eq:short_hand}
			\begin{align}
				\partial_t f_r + \nabla_x \cdot (vf_r) &= \nabla_v \cdot \left(S_r(f) f_r\right), \qquad
				\partial_t f_b + \nabla_x \cdot (vf_b) = \nabla_v \cdot \left( S_b(f) f_b \right), \\
				f_r(0,x,v) &= f_r(x,v),\qquad\qquad\qquad\qquad\, f_b(0,x,v) = f_b(x,v).
			\end{align}
		\end{subequations}
		All studies in the numerical section will be based on these equations.
	\end{example}
	The mean-field equations derived here are expected to be helpful for analytical investigations of stationary states and the formation of patterns in future work. Here, the main interest is to see whether or not the collision avoidance and pattern formation on the microscopic level translates to the mesoscopic system.

	\section{Boundary Conditions and Obstacles}\label{sec:boundaryObstacles}
	The purpose of modelling agent dynamics is to predict their behaviour in events of evacuation or in crowded sites. Therefore, it is reasonable to assume that the domain is bounded, and to be even more realistic, that there are obstacles in the domain. We begin with the discussion of boundary conditions.
	\subsection{Boundary conditions} 
\begin{figure}[htp]
		\begin{center}
		\includegraphics{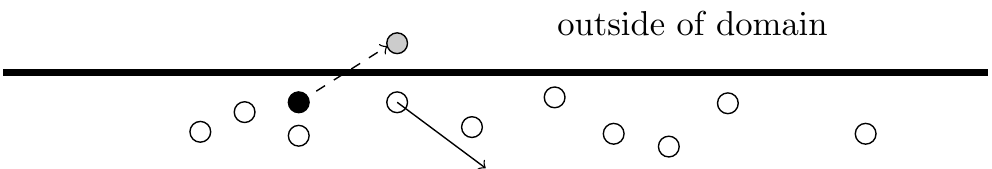}
		\caption{Illustration of the reflective boundary conditions. The black particle with velocity depicted with the dashed vector is about to leave the domain in the next time step. Due to the reflecting boundary conditions, it is projected into the domain and the y-component of its velocity is reflected (black vector).}
		\label{fig:boundaryConditions}
		\end{center}
	\end{figure}

	To model the fact that agents are in general not scared of walls, but usually keep some distance to a wall, we apply reflective boundary conditions. In fact, if an agent moves towards a boundary, he or she will be reflected. Thus, the velocity component which points towards the boundary will be transformed in a velocity pointing away from the boundary. Besides being realistic, these boundary conditions preserve mass, which is a desired property in this context, especially, in the mean-field limit. A sketch of the reflecting behaviour is given in Figure~\ref{fig:boundaryConditions}. 
	
	In order to study pattern formation, we assume to have periodic boundary conditions at the inflow and outflow of the domain. 
	\begin{example}
		For the simulations we consider a channel and a crossing. The domain and its periodic (yellow, green) and reflecting (black) boundaries are sketched in Figure~\ref{fig:boundarysetting}. The plot on the top refers to the channel setting, the one on the bottom to the crossing.
		\begin{figure}[htp]
			\begin{center}
			\includegraphics[scale=0.4]{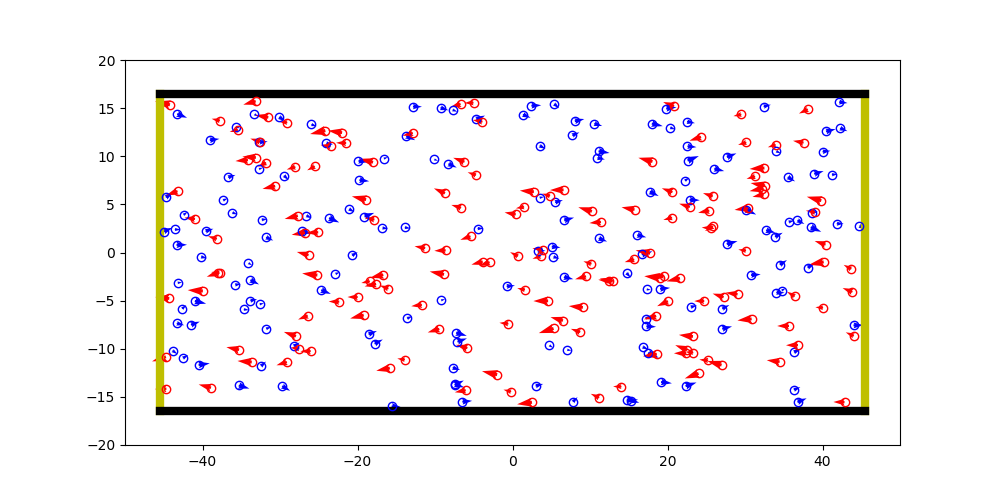}
			\includegraphics[scale=0.5]{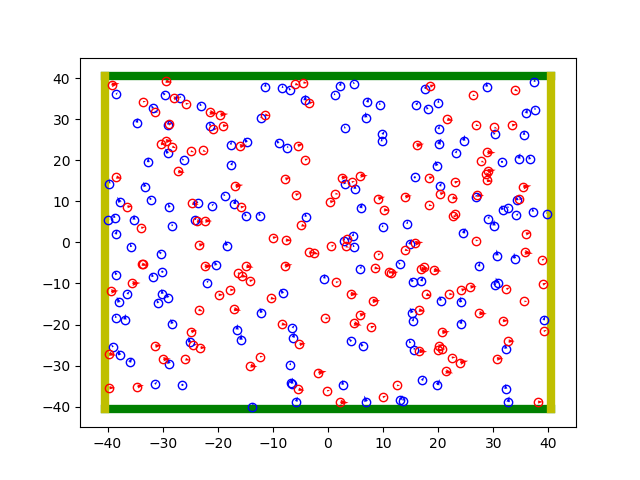}
			\caption{Illustration of the boundary setting. The green and yellow parts of the boundary refer to periodic boundary conditions. Walls are indicated by black lines. In case of the crossing the particles that leave the domain through a green boundary are flowing in at the other green boundary and analogous for the yellow boundary parts.}
			\label{fig:boundarysetting}
			\end{center}
		\end{figure}
	\end{example}
	
	\subsection{Obstacles}
\begin{figure}[htp]
		\begin{center}
		\includegraphics[scale=0.38]{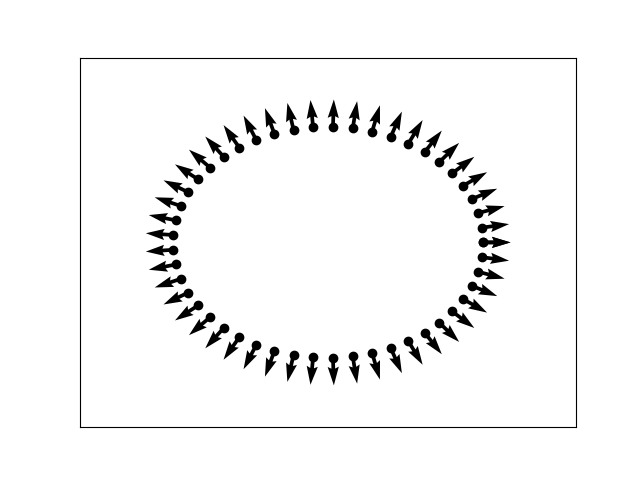}
		\includegraphics[scale=0.38]{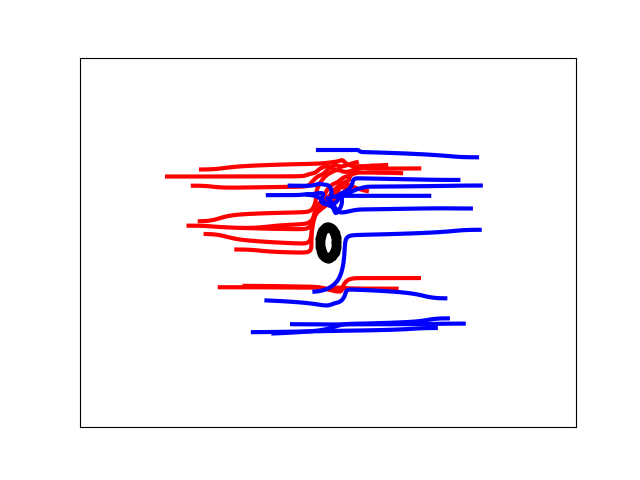}
		\caption{Left: Illustration of artificial agents with fixed position and artificial velocities modelling a circular obstacle. Right: Influence of a circular obstacle on the trajectories of the pedestrians.}
		\label{fig:illustraction obstacle}
		\end{center}
	\end{figure}

	The pairwise-interaction introduced above can be easily adapted to interaction with obstacles. Indeed, we define obstacles as artificial agents having a fixed position and an artificial velocity pointing outward of the obstacle. Positions and artificial velocity vectors of a circle shaped obstacle are illustrated in Figure~\ref{fig:illustraction obstacle} (left).
	The plot on the right shows the influence of this obstacle on trajectories of $20$ agents, $10$ red ones and $10$ blue ones.\\
	
	\section{Numerical schemes and implementation}\label{sec:numerics}
	In the following we discuss the numerical schemes used for the implementation on the microscopic and the mesoscopic scale. We begin with the particle implementation.
	\subsection{Particle Scheme}
	In order to solve the ODE system for the particles, we use a leap frog scheme combined with a splitting for the velocity update. In more detail, we compute
	\begin{align*}
		x_i^{k'} &= x_i^k + \frac\tau2 v_i^k, \qquad &&v_i^{k'} =   (v_i^k + \tau u_i)/(1+\tau) \\
		v_i^{k+1} &= v_i^{k'} + \tau M(v^{k'}) K(x^{k'}), \qquad &&x_i^{k+1} = x_i^{k'} + \frac\tau2 v_i^{k+1}.
	\end{align*}
	Here and in the following $\tau$ denotes the time step. Note that the part involving the desired velocity is solved implicitly and that the interaction is independent of the group membership of the two interacting particles. This allows for a straight forward vectorized implementation of the scheme.
	
	To include the periodic boundary conditions in the particle scheme, we use copies of the particles of interest. These artifical particles are representing the particles on the other end of the domain in the compuation of the interaction forces. Moreover, we restrict the domain of interaction in order to be consistent with the mean-field scheme. See Remark~\ref{rem:domain of interaction} below for more details.
	
		\subsection{Mean-field Scheme}
	For the mean-field simulation we employ a Strang splitting scheme. That means, a Semi-Lagrangian solver \cite{KlarReuterswaerd,Sonnendruecker} is used in the spatial domain and a Semi-Implicit Finite-Volume scheme is used in the velocity space.  The evolution is treated analogously for both groups. Using the short hand notation introduced in \eqref{eq:short_hand} and $i \in \{r,b\}$ we get
	\begin{subequations}\label{mf-splitting}
		\begin{align}\label{dis_1}
			\partial_t f_i^{k'} &= -\frac{1}{2} \nabla_v \cdot ( S_i(f^{k'})   f_i^{k'}), \qquad && f_i^{k'} = f_i(t),  \\
			\label{dis_2} \partial_t f_i^{k''} &= - v \cdot \nabla_x f_i^{k''}, \qquad && f_i^{k''} = f_i^{k'}(t + \tau),   \\
			\label{dis_3} \partial_t f_i^{k+1} &= -\frac{1}{2} \nabla_v  \cdot (S_i(f^{k'}) f_i^{k+1}), \qquad && f_i^{k+1} = f_i^{k''}(t+\tau) .
		\end{align}
	\end{subequations}
	The transport in \eqref{dis_2} is computed using a Semi-Lagrangian method. In fact, the computations are based on a fixed grid and  the mass is transported along characteristics. These characteristics curves are computed with the help of a second order Runge-Kutta scheme in a pre-processing step. In every transport step of the time loop, the initial point is a grid point. It is very unlikely to end again in a grid point, hence we need an interpolation to compute the new values on the grid points. For the interpolation we use a polynomial reconstruction with cubic Bezier curves which is of second order and has the nice property that it never leaves the convex hull of the control points. Additionally, we have to obey periodic or reflective boundary conditions in each step.
	
	For the transport in velocity space, i.e.~\eqref{dis_1} and \eqref{dis_3}, a second order finite volume scheme is employed. The advection is approximated by a Lax-Wendroff flux \cite{LeVeque, Quarteroni} and oscillations arising from non-smooth solutions are limited using a van-Leer method \cite{VanLeer}. In the velocity domain we do not need to fix boundary conditions, as we assume to never reach the boundary of the domain. More details on the schemes without boundary conditions can be found in \cite{Sheep1,CarrilloKlarRoth}. 
	\begin{remark}
		Note that in contrast to these references the interaction of the agents in this contribution is not only depending on their pairwise distance, but also on their velocities. This has to be taken care of in the finite volume scheme. Moreover, the dependence on the velocities increases the cost of computing the interaction forces significantly.
	\end{remark}
	
	In order to avoid cancellation, we make use of the structure $f = f_r + f_b$ and compute the interaction with each group separately. Hence, we do the procedure explained in \eqref{mf-splitting} first for the interaction with the red group and then for the blue group.
	
	\begin{remark}\label{rem:domain of interaction}
		Another important observation is the following: the interaction takes place in a very short range, thus to speed up the computation of the interaction force, we consider only points in the neighbourhood. Indeed, the domain of interaction at location ${i,j}$ is given by all the points $(k,l)$ with $k=i-2 \text{ to } i+2$ and $l=j-2\text{ to }j+2$. This speeds up the computation of the interaction forces significantly. 
		In order to make the simulations on the particle and the mean-field level comparable, we restrict the interaction domain in the particle scheme to a finite value as well.
	\end{remark}
	All codes for the particle simulations are written vectorized in Python 3.6, the mean-field simulations are implemented in Fortran 90 and their results are plotted with Matlab R2018a.
	
		\section{Numerical Results}\label{sec:simulations}
	In this section we discuss numerical results for the two settings sketched in Figure~\ref{fig:boundaryConditions} on the microscopic and on the mesoscopic level. We begin with the channel setting which has periodic boundary conditions on the left and the right side and reflective boundary conditions on the top and the bottom, see Figure~\ref{fig:boundarysetting} (top). Then we discuss the crossing with periodic boundary conditions at top and bottom and left and right, respectively. 
	
	The following simulation results are computed using interaction forces resulting from the Morse Potential \cite{Morse}
	$$ P(d) = Re^{-d/r} - Ae^{-d/a} $$
	with strength parameters $A=0, R = 500$ and range parameters  $r = 1.5, a = 1.5$. Further, $d$ refers to the distance of two interacting pedestrians $d = |x_i -x_j|.$ Moreover, we assume that all pedestrians are having the tendency to go to the right to avoid a collision, i.e., $\lambda = 0.25$.
	
	\begin{remark}
		Here $A = 0$ corresponds to the modeling assumption that the pedestrians have only repulsive influence on each other. Groups of friends or families are not considered here. A study of the influence of these on the simulations would be interesting as future work.
	\end{remark}
	
	\subsection{Pedestrians in a channel}
	The random initial positions are distributed uniformly in the domain $\Omega_x  = [-45,45] \times [-15,15]$. The initial velocities of the blue and the red particles are distributed in $$\Omega_{v_b} = [-0.3,-0.1] \times [-0.2,0.2]\quad \text{and}  \quad \Omega_{v_r} = [0.1,0.3] \times [-0.2,0.2],$$ respectively. Further parameters are
	$dt = 0.01,\; N_b \in  \{10, 150\},\; N_r \in \{10, 150\}$ at the particle level and $dt = 0.05,\; nx = 100,\; ny = 40,\; nv_1 = 20,\; nv_2 = 20$ for the mean-field simulation. Here $nx, ny, nv_1$ and $nv_2$ denote the number of discretization points for the spatial and velocity domain in $x$- and $y$-direction, respectively. Note that this choice implies that the red and the blue group have mass 0.5, as $f^N$ is by construction a probability measure. 
	
		\begin{figure}[htp]
		\begin{center}
			\includegraphics[scale=0.42]{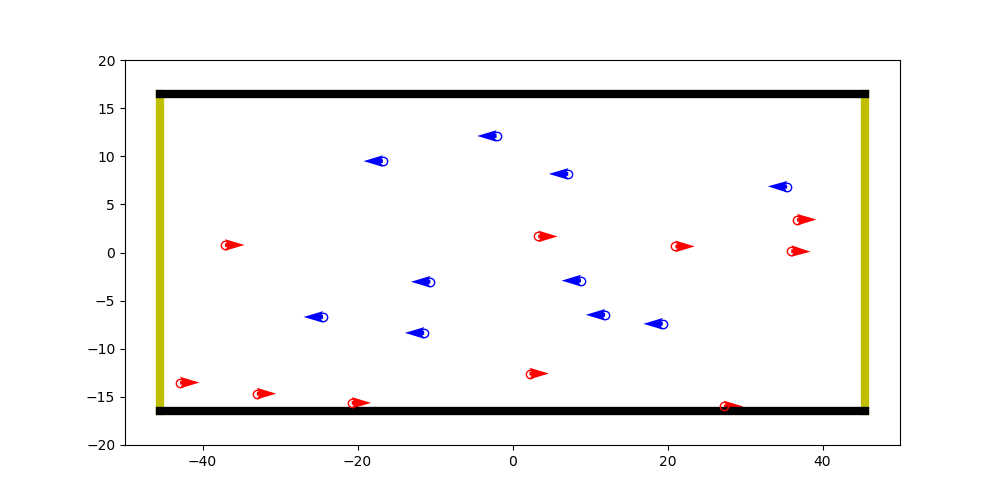}
			\caption{Lane formation in a channel - particle simulation with only few particles involved, we see a formation of multiple horizontal lanes as stationary state. 
			The parameters are $N_b=10=N_r, t=150.$ The arrows show the velocities of the particles.}
			\label{fig:lowdensity}
		\end{center}
	\end{figure} 
	
	In Figure~\ref{fig:lowdensity} we see positions of the groups with $N_r = 10 = N_b$ at time $t=600.$ Here, four horizontal lanes have formed and the interaction forces are small enough such that this is a stable configuration. For higher particle densities we expect a formation of only two lanes.
	\begin{remark}
		Numerical results suggest that there is a relation between the number of lanes formed and the parameter $r$ which scales the interaction radius of the particles. This seems to be in analogy to the investigation of the stationary number of opinions in opinion formation simulations, see \cite{MotschTadmor}. It would be very interesting to quantify the number of lanes for given interaction range $r.$
	\end{remark}

		\begin{figure}[htp]
		\begin{center}
			\includegraphics[scale=0.24]{channelStart.png}
			\includegraphics[scale=0.24]{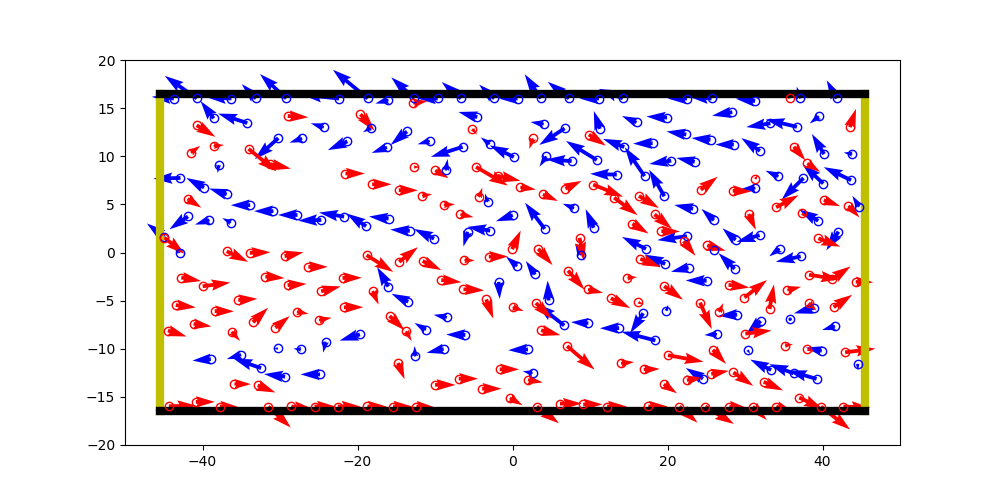}	\\
			\includegraphics[scale=0.24]{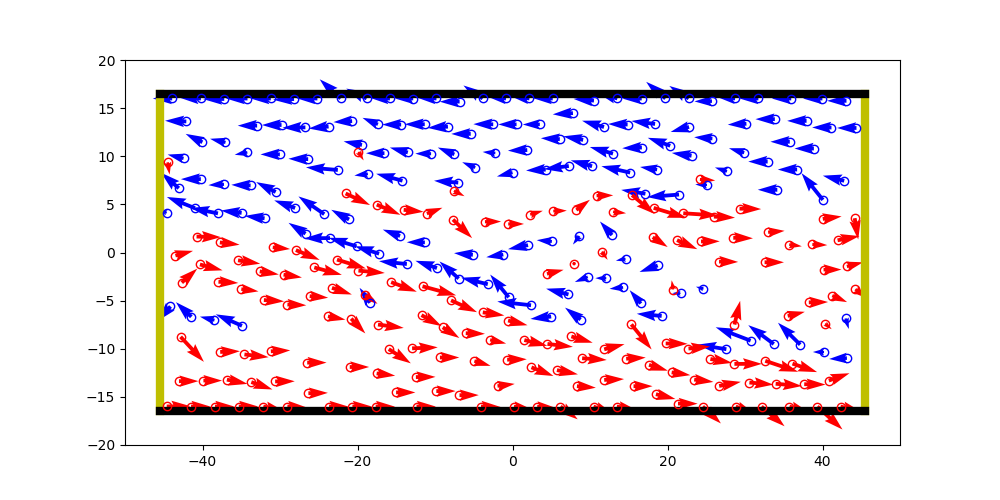}
			\includegraphics[scale=0.24]{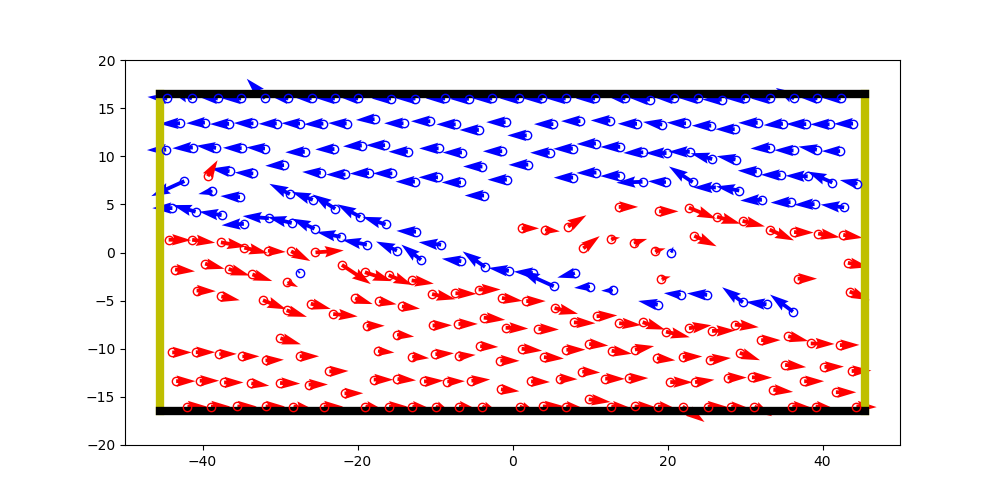}	\\
			\includegraphics[scale=0.24]{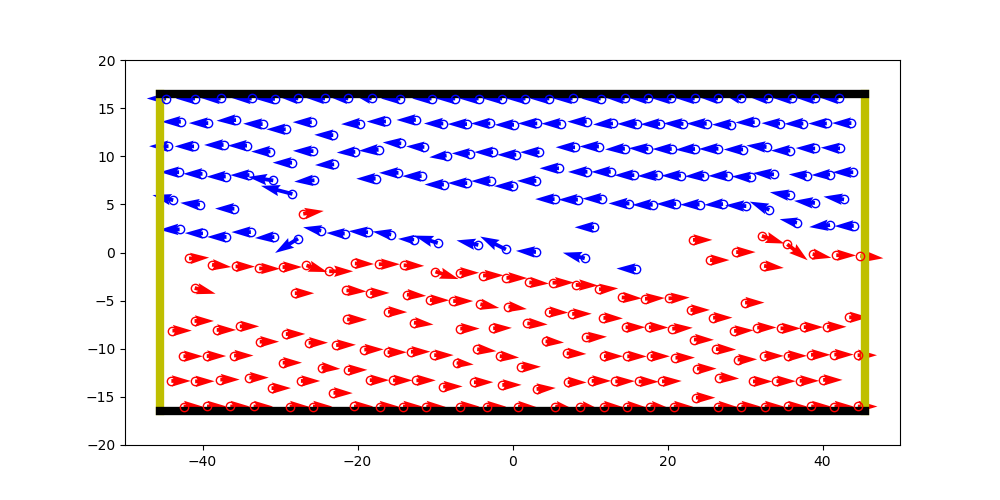}
			\includegraphics[scale=0.24]{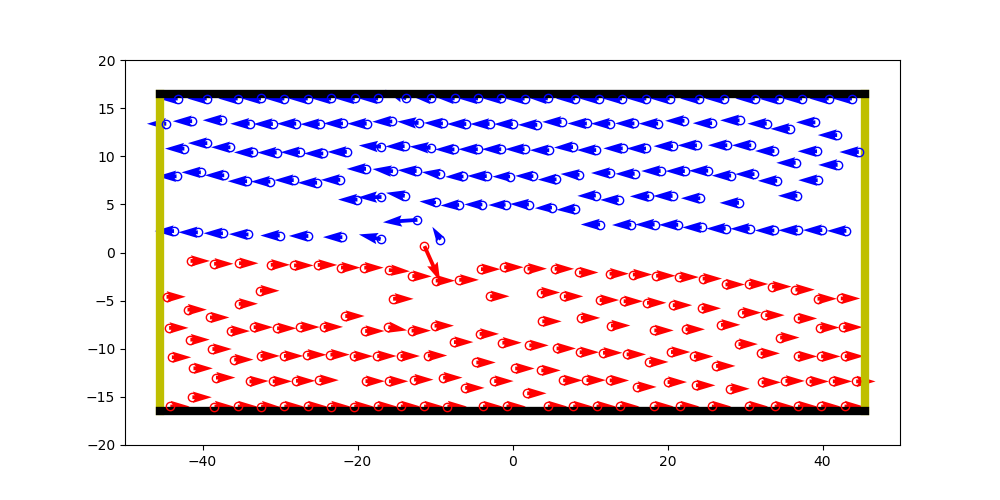}	
			\caption{Lane formation in a channel - particle simulation. The blue pedestrians go from right to left, the red from left to right, i.e. $u_b = (-0.2,0)^T$ and $u_r = (0.2,0)^T.$ The arrows show the velocities of the pedestrians. As $\lambda = 0.25,$ pedestrians prefer to step to the right to avoid a collision. The snapshots are made at the times $t=0,\; 50,\; 100,\; 150,\; 200,\; 250.$}
			\label{fig:laneformationParticles}
		\end{center}
	\end{figure} 
	
	Next, we investigate the channel setting with $N_b = 250 = N_r$ and the corresponding mean-field simulation. Figure~\ref{fig:laneformationParticles} shows the results of the particle simulation. Initially the particles are uniformly distributed over the whole domain as shown in the plot on the top-left. As time evolves we see the formation of diagonal stripes ($t=50,\;100$). The red particles are moving towards the bottom of the domain and the blue group is moving upwards. After some time a stationary distribution of two lanes, the red at the bottom of the domain and the blue at the top of the domain has formed ($t=250$).
	This is in contrast, to the behaviour shown in Figure~\ref{fig:lowdensity} suggesting that there is a correlation of interaction radius and the number of lanes. A detailed study on this is work in progress.
	
	\begin{remark}
		Heuristically, the formation of lanes can be explained by sequential occurrences of evasive behaviour. If we consider ourselves as one particle in the high density crowd, then soon we see someone approaching us. We need to adjust our velocity in order to avoid a collision. Then the next one approaches us, again, we avoid the collision. This behaviour terminates when we are at the top or bottom boundary of the domain, or, when we are surrounded by others walking in the same direction, hence, we feel no need to avoid any collision. The velocity vectors in Figure~\ref{fig:laneformationParticles} underline this thought experiment.
	\end{remark}	
	
	\begin{figure}[htp]
		\begin{center}
			\includegraphics[scale=0.17]{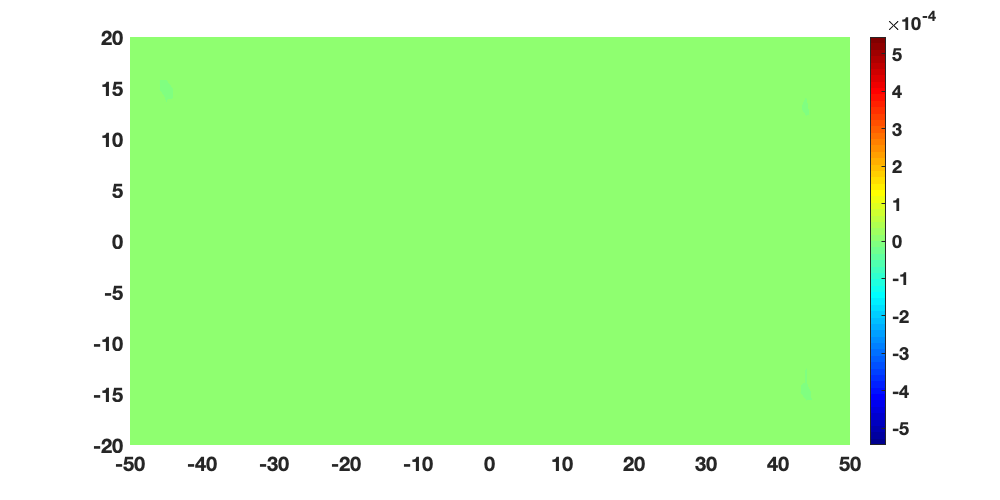}
			\includegraphics[scale=0.17]{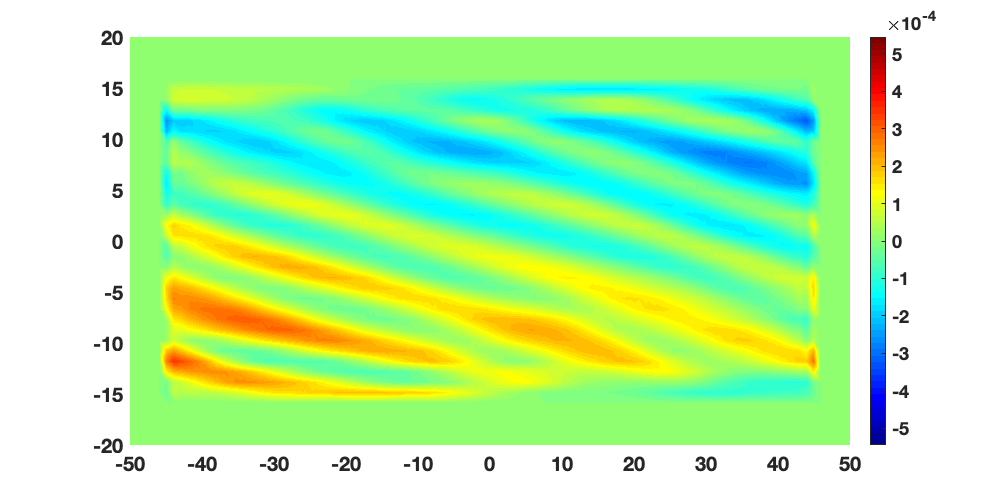} \\	
			\includegraphics[scale=0.17]{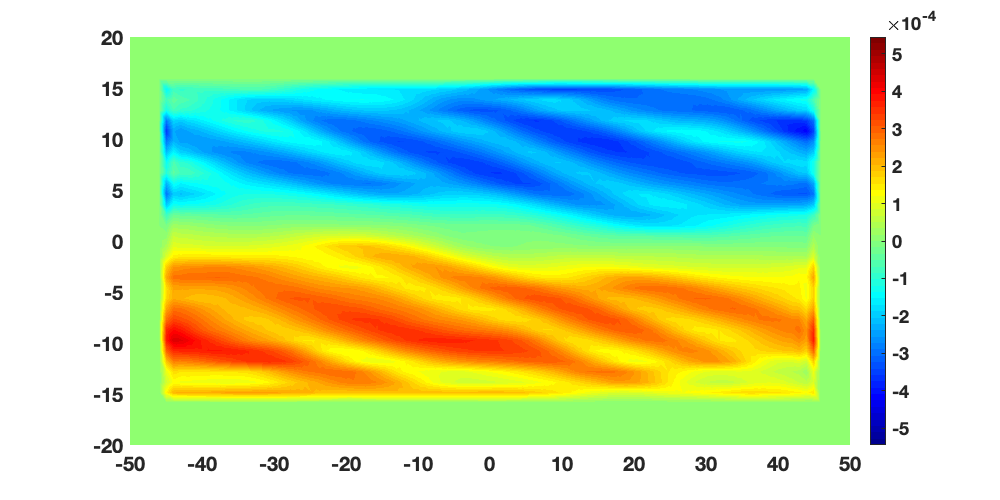}
			\includegraphics[scale=0.17]{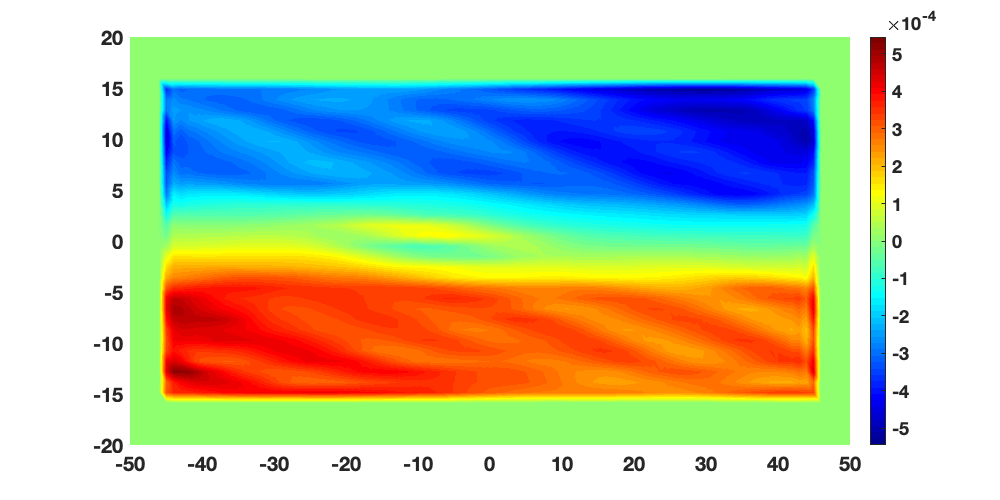} \\
			\includegraphics[scale=0.17]{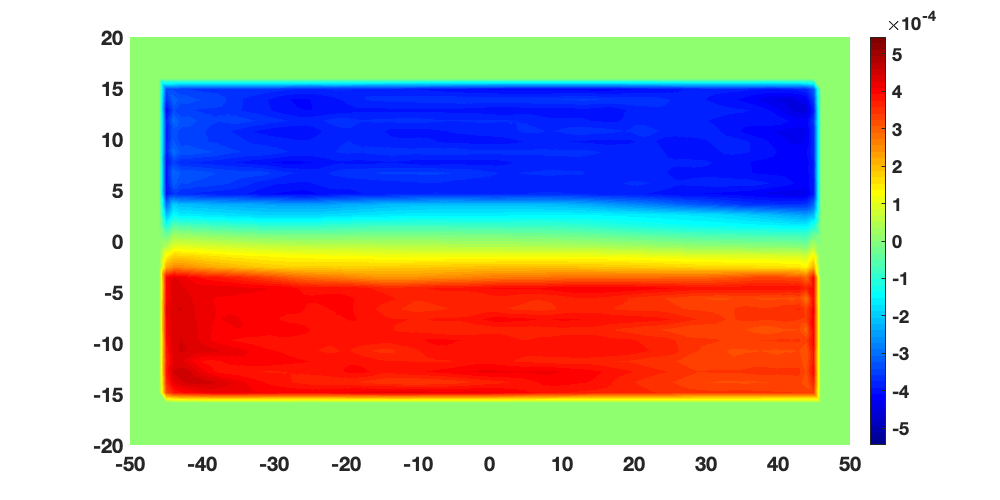}	
			\includegraphics[scale=0.17]{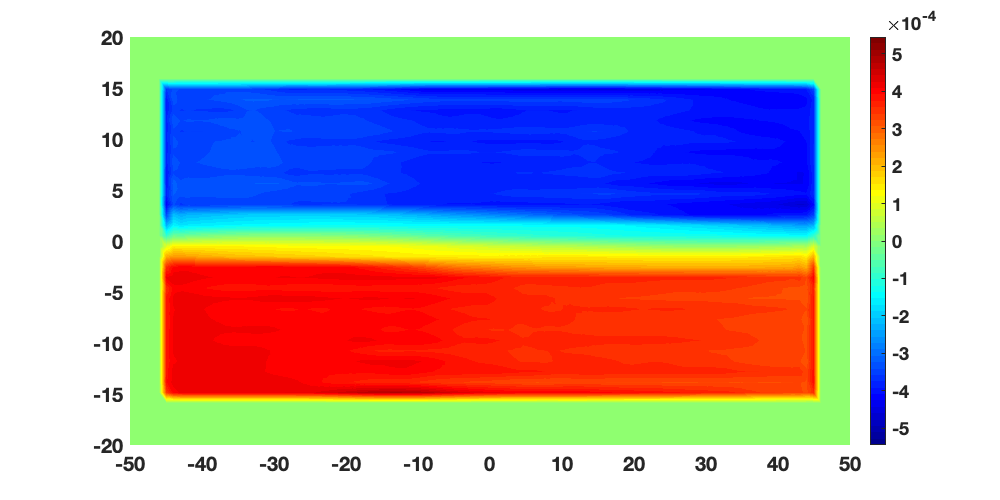}
			\caption{Lane formation in a channel - mean-field simulation. Initially the red and the blue group are uniformly distributed in the domain. The difference of the densities vanishes as shown in the plot on the top-left. Then we see the formation of diagonal stripes at time t=50 and finally a formation of lanes as time proceeds. The plots correspond to time times t=0, 250, 500, 750, 1000, 1250.}
			\label{fig:laneformationmeanfield}
		\end{center}
	\end{figure} 
	
	The same observations are made for the simulation on the mean-field level. To discuss the results, we introduce
	\begin{equation}\label{quantities}
		\rho_i(t,x) := \int f_i(t,x,v)\;  \mathrm{d}v\quad \text{and}\quad \phi_i(t,v) := \int f_i(t,x,v) \; \mathrm{d}x, \qquad i=\{r,b\}. 
	\end{equation}	
	Figure~\ref{fig:laneformationmeanfield} shows the quantity $\rho_r(t,\cdot) - \rho(t, \cdot)$ at several time instances $t$. Initially the quantity is zero, as both groups are uniformly distributed in the spatial domain. As the groups start to cluster and to form lanes, we see the density of the blue group depicted in blue and the density of the red group in red. As expected from the formal derivation above, the results match the results of the particles level.

\begin{figure}[htp]
		\begin{center}
		\includegraphics[scale=0.2]{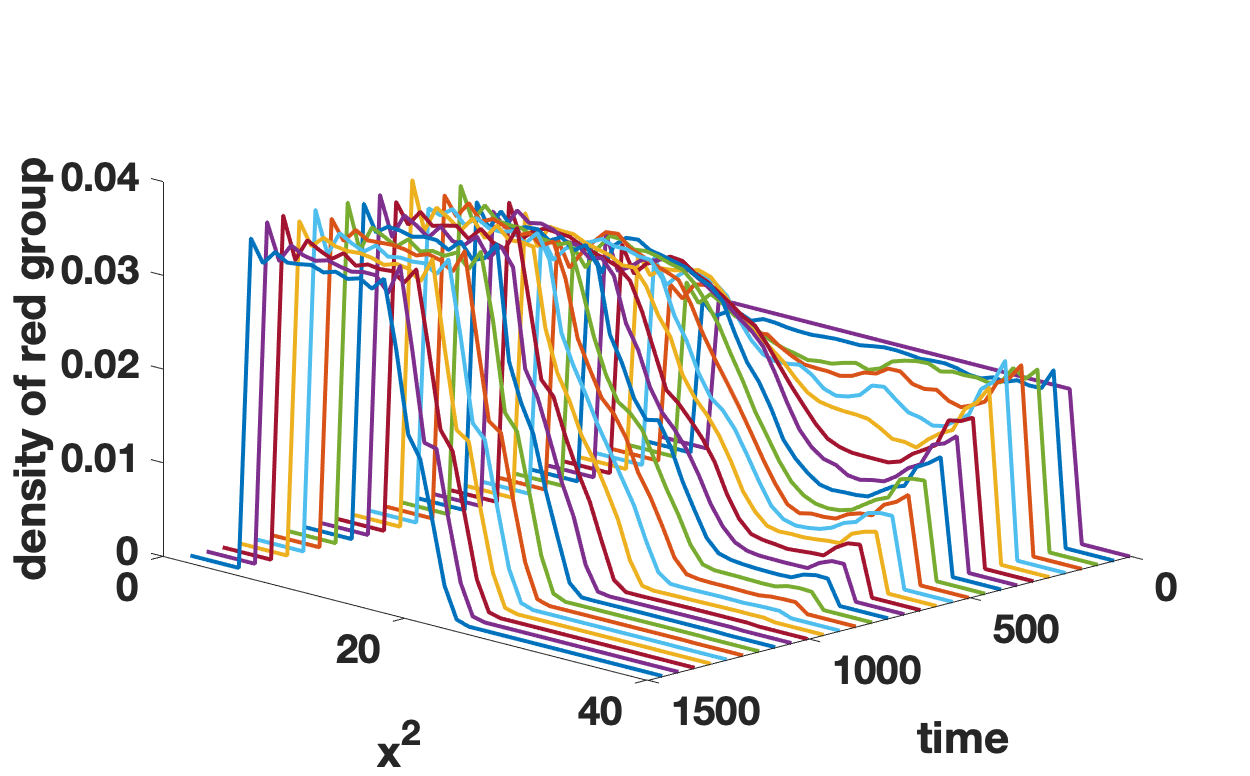}
		\hfill
		\includegraphics[scale=0.2]{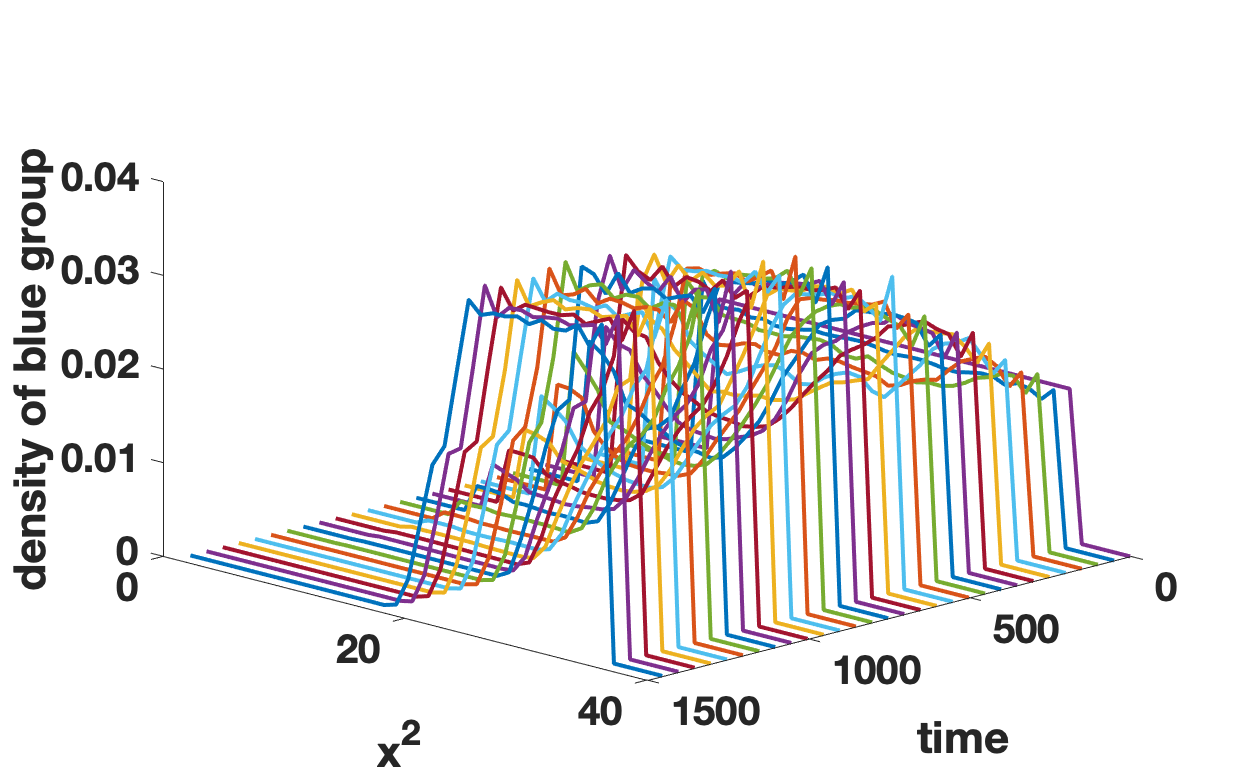}	
		\caption{Density distribution in a channel - mean-field simulation. The densities $\rho_i^2, i \in \{r, b\}$ are integrated w.r.t. $x^1$ in order to extract the information along the $x^2$-axis.}
		\label{fig:densitydistribution}	
		\end{center}
	\end{figure}

\begin{figure}[htp]
		\begin{center}
		\includegraphics[scale=0.3]{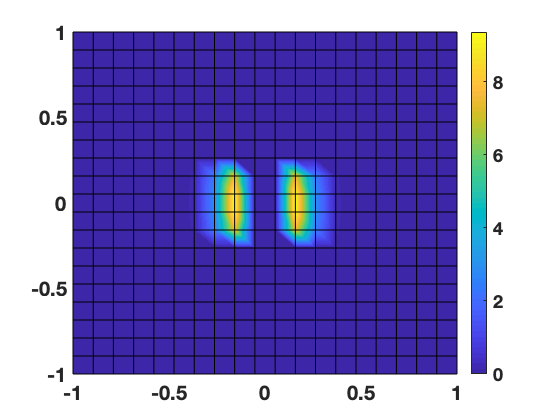}
		\includegraphics[scale=0.3]{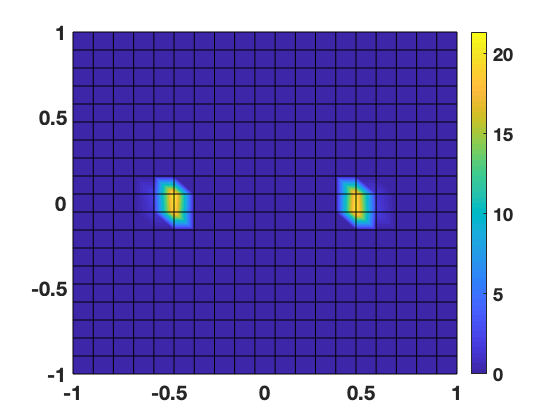}	
		\includegraphics[scale=0.3]{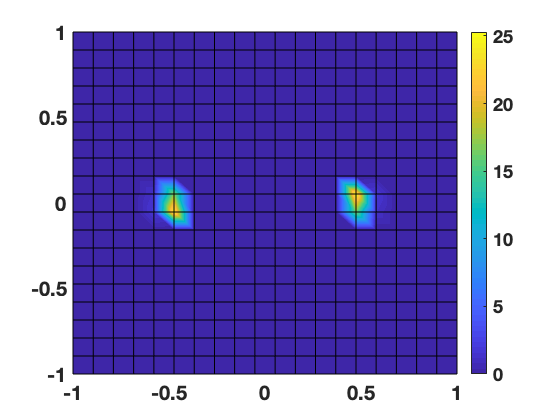}
		\caption{Lane formation in a channel - mean-field simulation. Initially the red and the blue group are uniformly distributed in the velocity domain $\Omega_{v_r}$ and $\Omega_{v_b}$, respectively, as shown in the plot on the top-left. Then the relaxation towards the desired velocities $u_r = (0.2,0)^T$ and $u_b = (-0.2,0)^T$ starts, see plot on the top-right. Only small changes in the $y$-coordinates lead the crowd to the lane configuration as can be seen when comparing the top-left and the bottom plot of the velocity density. The time instances are $t=0,\; 250,\; 1500.$}
		\label{fig:velocitychannel}
		\end{center}
	\end{figure}

	The evolution of $$\rho_i^2(t,x^2) = \int f_i(t,x^1,x^2,v_1,v_2)\; \mathrm{d}(v_1,v_2)\; \mathrm{d}x^1\qquad i\in \{r,b\}$$ is illustrated in Figure~\ref{fig:densitydistribution}. Note that the quantities are averaged with respect to the $x^1$-component as well in order to emphasize the evolution in the $x^2$-component.
	Finally, we study the evolution of the velocity density $\phi_r(t,\cdot) + \phi_b(t,\cdot)$ in Figure~\ref{fig:velocitychannel}. We see the initial distribution at the top-left plot. The relaxation towards the desired velocities is achieved after a very short time. Then only small changes in the $x^2$-component lead to the formation of the lanes.
	
	\begin{remark}
		For negative anisotropy parameter, i.e. $\lambda=-0.25$ the roles of the to groups change. In this case, the red lane forms at the top of the domain and the blue lane at the bottom of the domain.
	\end{remark}	
	
		\subsection{Pedestrians at a crossing}
\begin{figure}[htp]
		\begin{center}
		\includegraphics[scale=0.35]{crossStart.png}
		\includegraphics[scale=0.35]{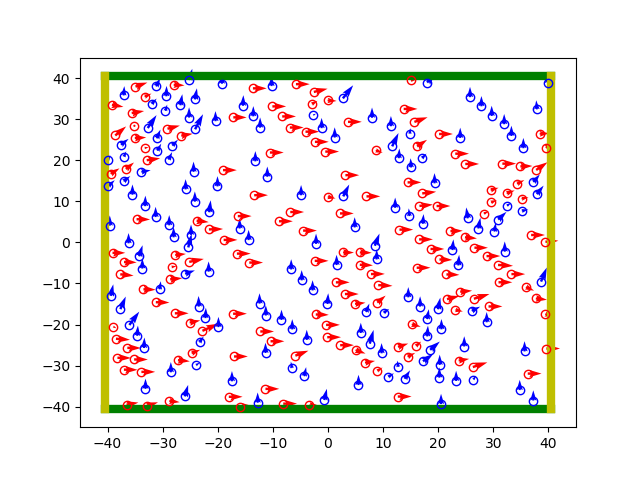}	
		\includegraphics[scale=0.35]{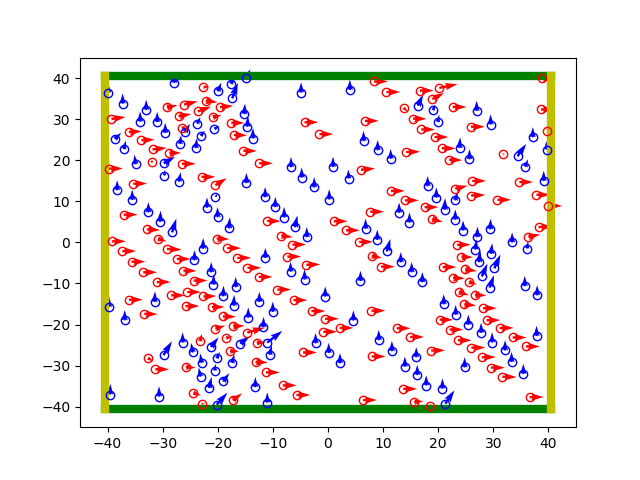}
		\includegraphics[scale=0.35]{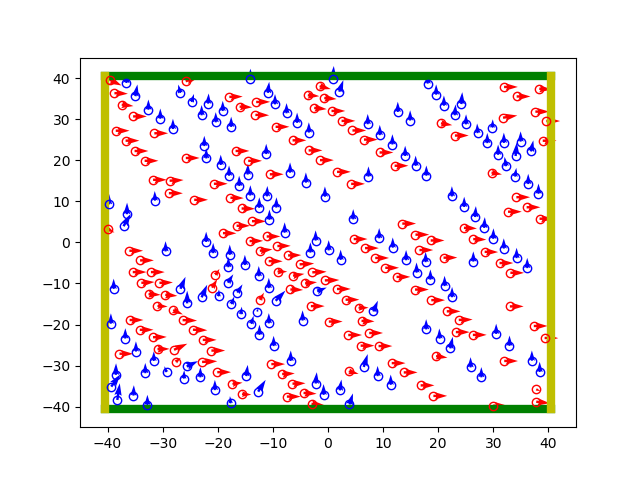}
		\includegraphics[scale=0.35]{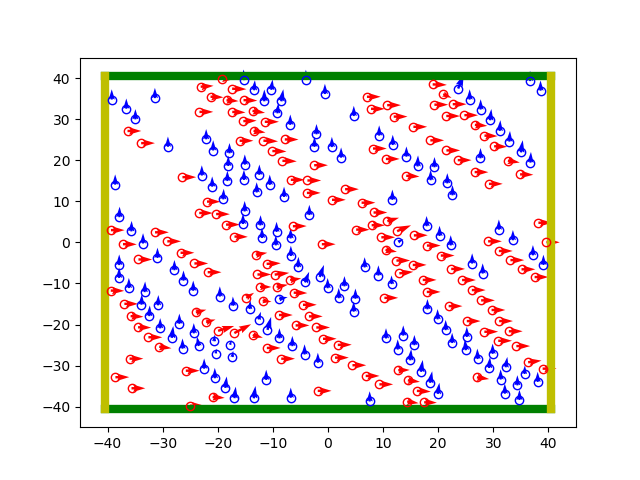}	
		\includegraphics[scale=0.35]{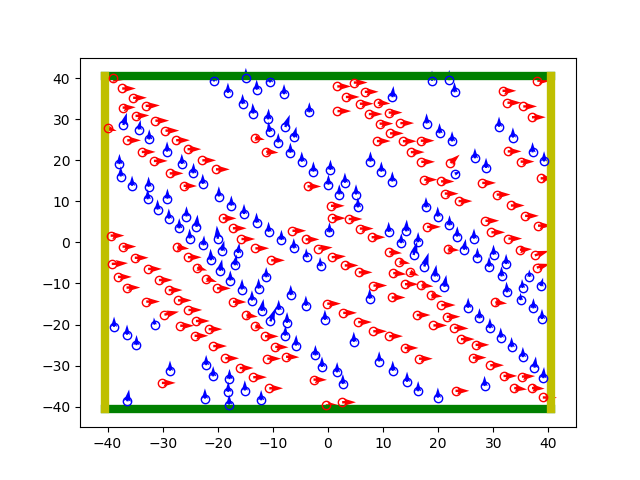}
		\caption{Travelling waves at a crossing - particle simulation. The blue go from bottom to top, the red from left to right, i.e. $u_b = (0,0.2)$ and $u_r = (0.2,0).$ The arrows show the velocities of the pedestrians. As $\lambda = 0.25,$ pedestrians prefer to step to the right to avoid a collision. From top-left to bottom-right: $t = 0$, $t = 50$, $t = 100$, $t = 150$, $t = 200$ and $t=250$. The simulation was done with $N_b=150=N_r.$ At $t=250$ the pattern allows most pedestrians to walk with their desired velocities.}
		\label{fig:travellingwavesParticles}
		\end{center}
	\end{figure}

	For the simulation of the crossing the initial positions and velocities are randomly distributed in $\Omega_{x} = [-40,40] \times [-40,40]$ and $\Omega_{v} = [-0.1,0.1] \times [-0.1,0.1]$, respectively. The desired velocity for the red group is $u_r = (0.2,0)^T$ and for the blue group it is $u_b = (0,0.2)^T.$ All other parameters are as in the previous case, except for $nx = 40 = ny.$ Moreover, boundary conditions are periodic from left to right and from top to bottom as indicated by the green and yellow lines in Figure~\ref{fig:travellingwavesParticles}.
	
	The figure shows corresponding simulation results for the particle scheme. Initially the particles are distributed uniformly, after some time we see patterns forming. This process ends in stable travelling waves as shown in the plot at the bottom-right ($t=250$).
	
	\begin{figure}[htp]
		\begin{center}
		\includegraphics[scale=0.3]{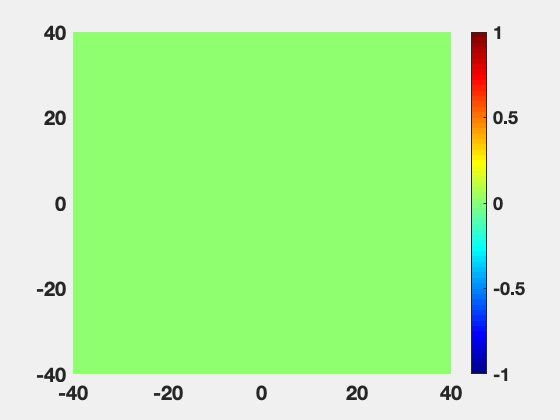}
		\includegraphics[scale=0.3]{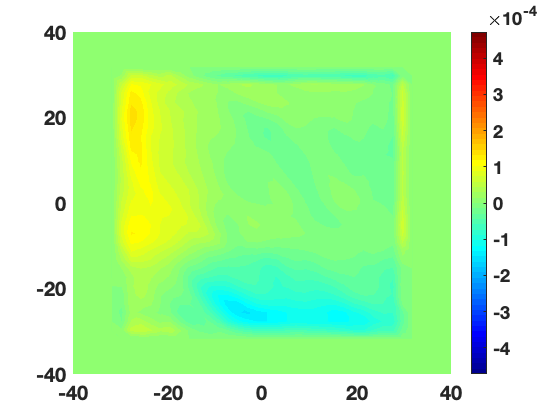}	
		\includegraphics[scale=0.3]{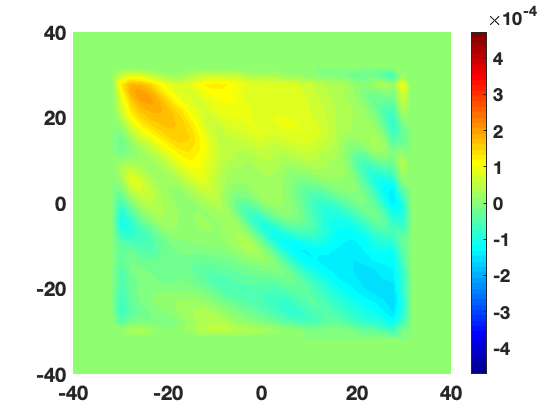}
		\includegraphics[scale=0.3]{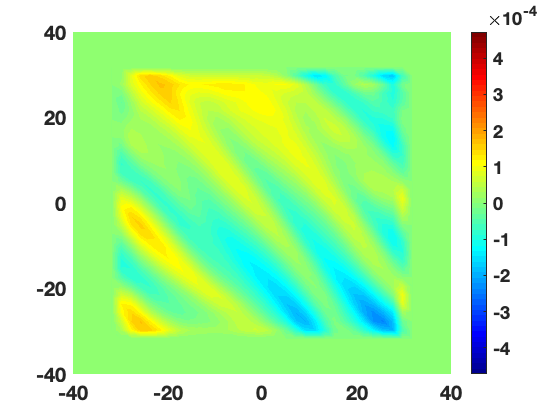}
		\includegraphics[scale=0.3]{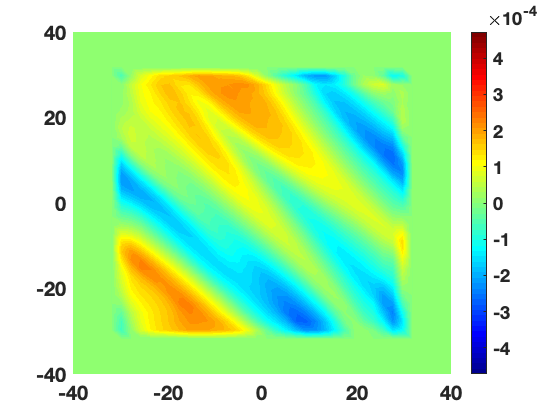}	
		\includegraphics[scale=0.3]{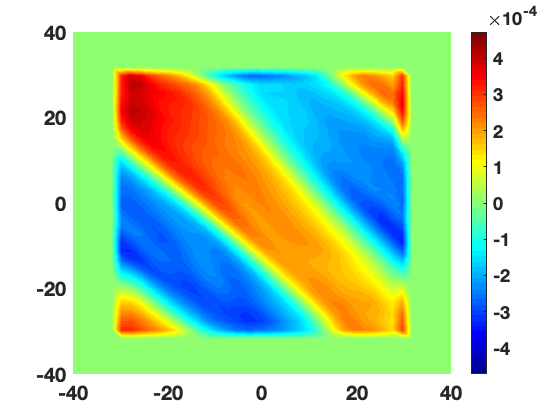}
		\caption{Travelling waves at a crossing - mean-field simulation. Initially the red and the blue group are uniformly distributed in the spatial domain. The difference of the densities vanishes as shown in the plot on the top-left. At time t=50 the groups start to separate, see the plot in the top-right. Afterwards the line pattern is forming (t=150, 375 and 750). Finally, we see the stationary configuration of travelling waves at time t=1500 in the bottom-right plot.}
		\label{fig:travellingwavesmeanfield}
		\end{center}
	\end{figure} 
	
\begin{figure}[htp]
		\begin{center}
		\includegraphics[scale=0.3]{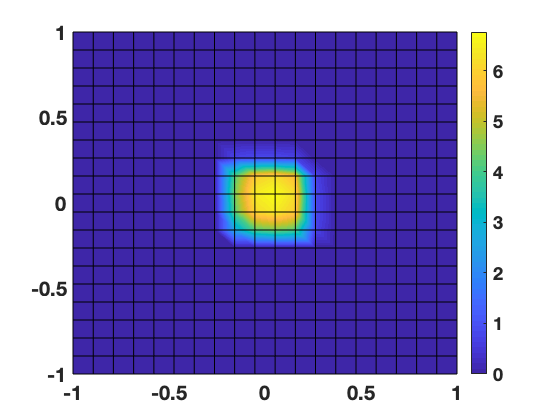}
		\includegraphics[scale=0.3]{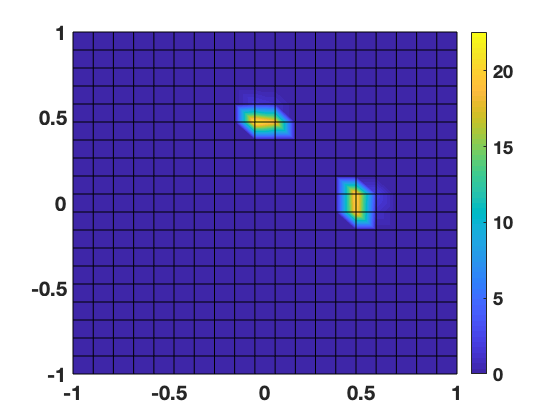}	
		\includegraphics[scale=0.3]{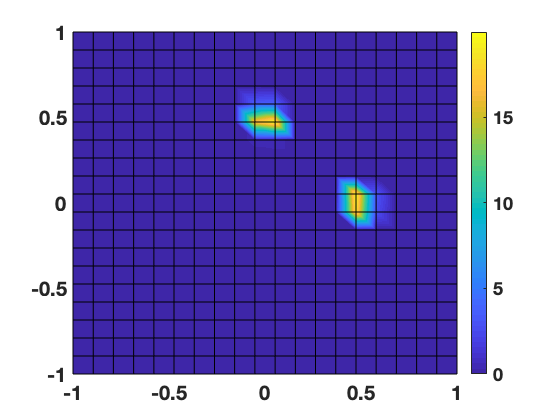}
		\caption{Travelling waves at a crossing - mean-field simulation. Initially the red and the blue group are uniformly distributed in the velocity domain $\Omega_v$ as shown in the plot on the top-left. Then the relaxation towards the desired velocities $u_r = (0.2,0)^T$ and $u_b = (0, 0.2)^T$ starts, see plot on the top-right. Only small changes in the $x$- and $y$-coordinate lead the crowd to the travelling waves configuration as can be seen when comparing the top-left and the bottom plot of the velocity density. The time instances are $t=0, 250, 1500.$}
		\label{fig:velocitycross}
		\end{center}
	\end{figure}

	Simulation results for the mean-field setting are depicted in Figure~\ref{fig:travellingwavesmeanfield} and Figure~\ref{fig:velocitycross}. The plots show the spatial and velocity densities of the red and the blue group, which are defined as in the channel setting, see \eqref{quantities}.
	In Figure~\ref{fig:travellingwavesmeanfield} we see a plot of the difference $\rho_r(t,x) - \rho_b(t,x)$ at various time instances $t.$ Hence, in locations where the blue group is clustered, we see blue colour and in locations with a majority of the red group we see red colour. At the beginning both groups are uniformly distributed over the spatial domain, the difference is therefore zero, as shown in the top-left plot. As time evolves we see the formation of groups which ends in a stable configuration of travelling waves as can be seen in the bottom-right plot of Figure~\ref{fig:travellingwavesmeanfield}. Thanks to the travelling wave pattern, almost every pedestrian walks with his/her desired velocity. Note that these patters are also reported in many other simulations and experiments with crossing flows, see for example \cite{Appert-Rolland,Ranetbauer1,Cividini,Ranetbauer2}.
	
	The velocity density $\phi_r(t,\cdot) + \phi_b(t,\cdot)$ is shown in Figure~\ref{fig:velocitycross} for difference time instances $t$.  Initially the red and the blue group are uniformly distributed in the velocity domain $\Omega_v$ as shown in the plot on the top-left. Then the relaxation towards the desired velocities $u_r = (0.2,0)^T$ and $u_b = (0, 0.2)^T$ starts, see plot on the top-right. Only small changes in the $v^1$- and $v^2$-coordinate lead the crowd to the travelling waves configuration as can be seen when comparing the top-left and the bottom plot of the velocity density. 
	
	\section{Conclusion and Outlook}\label{sec:conclusionOutlook}
	We discussed an anisotropic interaction model that allows for collision avoidance. The model is based on general swarm models and uses a rotation matrix to allow for a smooth evasion behaviour. First, we looked at a binary interaction of two agents, then we stated the particle system for $N$ agents and derived formally a mesoscopic formulation. We illustrated how obstacles can be included in the model conveniently. The behaviour of the model is studied numerically with the help of two scenarios. The first one considers a channel with two groups of pedestrians, one group walking from right to left and the other group from left to right. We see the typical lane formation observed in many other works and in reality. The second scenario investigates the behaviour of pedestrians at a crossing. Again, we see the typical travelling wave formations on both, the particle and the mean-field level. 
	
	In a next step, the model shall be analysed: a study of well-posedness in the ODE and PDE level is planned. Further, a rigorous proof of the relation of the particle and the mean-field system has $N\rightarrow \infty$ is to be achieved. Finally, it is interesting to find an analytical way to show the pattern formation as $T\to \infty$. Based on this work, an optimization framework corresponding to an evacuation scenario and a parameter identification based on data is planned for further studies.
	
	\section*{Acknowledgements} The author is gratefully acknowledging the motivating feedback given by Martin Burger and Marie-Therese Wolfram on the first sketch of the idea leading to this article.

\end{document}